\documentclass[11pt]{article}

\usepackage[top = 2cm, bottom = 2cm, left = 1.8cm, right = 1.8cm, paperheight=24cm,paperwidth=17cm]{geometry}
\usepackage{times}
\usepackage{epsfig} 
\usepackage{amsmath}
\usepackage{graphics}
\usepackage{float}
\usepackage{xcolor}
\usepackage[normalem]{ulem}

\usepackage[normalsize, it, bf]{caption}
\renewcommand \thesection{\arabic{section}}
\usepackage[utf8x]{inputenc}
\usepackage[greek,english]{babel}

\makeatletter
\renewcommand{\@cite}[1]{#1}

\renewcommand\thesubsection   {\thesection\@arabic\c@subsection}

\newcommand\sectionn{\@startsection {section}{1}{\z@}%
	{-3.5ex \@plus -1ex \@minus -.2ex}%
	{2.0ex \@plus.2ex}%
	{\centering\normalfont\large\bfseries}}
\def\list#1#2{\ifnum \@listdepth >5\relax \@toodeep
	\else \global\advance\@listdepth\@ne \fi
	\rightmargin \z@ \listparindent\z@ \itemindent\z@
	\csname @list\romannumeral\the\@listdepth\endcsname
	\def\@itemlabel{#1}\let\makelabel\@mklab \@nmbrlistfalse #2\relax
	\parskip 0pt\itemsep 0pt\parsep 0pt \topsep 0pt%
	\@trivlist
	\parindent\listparindent
	\advance\linewidth -\rightmargin \advance\linewidth -\leftmargin
	\advance\@totalleftmargin \leftmargin
	\parshape \@ne \@totalleftmargin \linewidth
	\ignorespaces}

\def\@startsection#1#2#3#4#5#6{%
	\if@noskipsec \leavevmode \fi
	\removelastskip
	
	\par
	
	\@tempskipa #4\relax
	\@afterindenttrue
	\ifdim \@tempskipa <\z@
	\@tempskipa -\@tempskipa \hspace{-2cm} \@afterindentfalse
	\fi
	
	\if@nobreak
	\everypar{}%
	\else
	\addpenalty\@secpenalty\addvspace\@tempskipa
	\fi
	\@ifstar
	{\@ssect{#3}{#4}{#5}{#6}}%
	{\@dblarg{\@sect{#1}{#2}{#3}{#4}{#5}{#6}}}}

\def\@sect#1#2#3#4#5#6[#7]#8{%
	\ifnum #2>\c@secnumdepth
	\let\@svsec\@empty
	\else
	\refstepcounter{#1}%
	\protected@edef\@svsec{\@seccntformat{#1}\relax}%
	\fi
	\@tempskipa #5\relax
	\ifdim \@tempskipa>\z@
	\begingroup
	#6{%
		\@hangfrom{\hskip #3\relax\@svsec}%
		\interlinepenalty \@M #8\@@par}%
	\endgroup
	\csname #1mark\endcsname{#7}%
	\addcontentsline{toc}{#1}{%
		\ifnum #2>\c@secnumdepth \else
		\protect\numberline{\csname the#1\endcsname}%
		\fi
		#7}%
	\else
	\def\@svsechd{%
		#6{\hskip #3\relax
			\@svsec #8}%
		\csname #1mark\endcsname{#7}%
		\addcontentsline{toc}{#1}{%
			\ifnum #2>\c@secnumdepth \else
			\protect\numberline{\csname the#1\endcsname}%
			\fi
			#7}}%
	\fi
	\@xsect{#5}}

\def\section{\@startsection{section}{1}{\z@}{12pt}{6pt}{\bf}}

\def\subsection{\@startsection {subsection}{2}{\z@}{12pt}{6pt}{\bf}}

\newtheorem{thm}{Theorem}[section]

\usepackage{abstract}
\usepackage{etoolbox}
\patchcmd{\abstract}{-.5em}{0em}{}{}
\AtBeginDocument{\parindent = 1.5em}

\providecommand{\keywords}[1]
{
	{\raggedright
	} \textit{Keywords:} #1}

\usepackage{titlesec}
\titleformat{\section}
{\normalfont\Large\bfseries}{\thesection}{3pt}{}
\usepackage{titlesec}
\titleformat{\subsection}
{\normalfont\Large\bfseries}{\thesubsection}{3pt}{}

\begin{document}
        
\begin{centering}
   \begin{minipage}{17cm}
       \selectlanguage{english}
   \end{minipage}
\end{centering}
\begin{centering}
   \vspace*{11pt} \vspace*{11pt} \vspace*{11pt} 
   \vspace{-.31cm}

   \Large\textbf{REVENUE AND SOCIAL WELFARE MAXIMIZATION IN SERVICE SYSTEMS: INDEPENDENT VERSUS GROUPED CUSTOMERS\normalsize\\} \vspace{11pt} \vspace{12pt}
   \textbf{\textit{A. Burnetas$^1$, A. Economou$^1$}}\\
   \small
   $^1$ Department of Mathematics, National and Kapodistrian University of Athens, Greece\\
   \{aburnetas, aeconom\}@math.uoa.gr
              
   \normalsize
\end{centering}
\vspace{11pt}
\normalsize
\thispagestyle{empty}

\begin{abstract}
\vspace*{-6pt}
   Revenue maximization in service systems with strategic customers is a central topic in the rational queueing literature. We study the join-or-balk problem in an unobservable M/M/1 queue under two behavioral regimes. Under independent behavior, customers choose individually whether to join; under grouped behavior, a customer leader chooses the joining fraction to maximize aggregate customer surplus. We show that, at any given fee, grouped behavior weakly reduces participation relative to independent behavior. The same ordering holds for the revenue-maximizing joining fractions and fees. Under independent behavior, the provider's revenue-maximizing joining fraction coincides with the social optimum and the provider extracts all customer surplus. Under grouped behavior, revenue maximization generally leads to lower participation and lower provider revenue, while leaving customers with positive surplus.
 \end{abstract}
\vspace*{6pt}

\setlength{\parindent}{13pt}

        
\section{INTRODUCTION}\label{introduction}
        
The study of service systems from a microeconomics viewpoint goes back at least to the pioneering papers of Naor (1969) and Edelson and Hildebrand (1975) who studied the join-or-balk dilemma for the customers in the M/M/1 queue, taking into account the strategic nature of the customers and of the service provider, who interact with the objective of maximizing their own utilities. These works recognized that the study of strategic interaction among customers and firms in the framework of service systems should take into account the presence of waiting costs in addition to the usual service values and costs that appear in other economic systems. Therefore, a credible economic study of a service system requires the amalgamation of ideas from classical Microeconomics with Operations Research, in particular Queueing Theory. This thread of the literature has evolved into a mature subarea of Queueing Theory that is usually referred to as Rational Queueing. Hassin and Haviv (2003) provided a nice overview of the basic methodology and the early results in this area. The monographs by Stidham (2009) and Hassin (2016) contain further material and overviews about models and methodologies in this subfield of Queueing Theory. A short introduction with an emphasis on the role of information that is provided to the customers has been presented in Economou (2020).

Service systems with strategic customers are commonly modeled using the join/balk framework, where arriving customers decide whether to enter a queue based on expected benefits and waiting costs. A standard assumption in this literature is that customers act independently, each optimizing their own expected utility while taking system congestion as given. This leads to an equilibrium joining behavior that determines the effective arrival rate and, in turn, system performance.

While analytically convenient, the assumption of independent decision making may not always be appropriate. In many settings, customers may coordinate their actions, either explicitly or implicitly, and make participation decisions collectively. Such coordination may arise through centralized control, platform-mediated decisions, or shared information and incentives. In these cases, customers effectively internalize the congestion they impose on one another, leading to fundamentally different behavior.

In this paper, we examine the join–or–balk problem in an unobservable single-server Markovian queue under two alternative behavioral regimes: independent customers and grouped (cooperating, coordinated) customers. In the former case, customers act individually and the system is characterized by a standard equilibrium. In the latter, a coordinating agent determines the participation level so as to maximize the collective benefit of customers. This leads to a distinct equilibrium notion, which can be interpreted as a socially coordinated outcome among customers.

Our main results provide an explicit ordering of the two regimes. At any fixed entrance fee, the joining fraction under grouped behavior is no greater than under independent behavior. The revenue-maximizing joining fraction, entrance fee, and provider revenue are also weakly lower under grouped behavior. Moreover, under independent behavior the provider's revenue-maximizing joining fraction coincides with the social optimum, whereas under grouped behavior revenue maximization generally induces too little participation relative to the social optimum but leaves customers with positive surplus.

The distinction between independent and coordinated decision making is closely related to the classical contrast between user equilibrium and system-optimal solutions in congestion models, where decentralized behavior fails to internalize congestion externalities. In a similar spirit, coordinated customer behavior in queueing systems can be interpreted as a collective decision that internalizes waiting costs.

Such coordination may arise through centralized control, platform-mediated decisions, or other mechanisms. For example, in Benioudakis et al. (2023), an aggregator platform essentially acts as a collective customer who responds to a service provider's service fee by setting a second price to his own customers. In a different setting, Bountali et al. (2022) consider an unobservable M/M/1 queue with strategic customers who arrive in batches of random size. Cooperating behavior is manifested in this model in the case where each arriving batch makes a collective decision  to join the system or balk, either as a whole or in part. 

From an economic perspective, coordinated customer behavior can be interpreted as a form of joint decision making that aligns individual incentives with collective outcomes, analogous to collusive arrangements in industrial organization where agents maximize joint payoffs rather than acting non-cooperatively (see e.g. Tirole (1988) Chapter 6).

The remainder of the paper is organized as follows. In section \ref{Section-Model} we present the model and its operational and economic parameters and state the pertinent optimization problems. In section \ref{Section-Anal-Indep-Cust} we summarize the results that have been reported in the literature for the case of independent customers. In section \ref{Section-Anal-Colluded-Cust} we consider the same problems for the case of grouped customers and develop the corresponding results: form of equilibrium strategies, revenue, customers' surplus and social welfare. In section \ref{Section-Comparison}, we compare the two cases and, in section \ref{Section-Computational}, we provide a number of findings from various numerical experiments. Conclusions and directions for further research are presented in section \ref{Section-Conclusions}.

\section{THE MODEL}\label{Section-Model}

We consider the service system of the standard single-server Markovian queue, also known as an M/M/1 queue. In this model, customers arrive according to a Poisson process with rate $\lambda$ and service times are exponentially distributed with rate $\mu$. There is a single server who serves customers one-by-one according to their order of arrivals (FCFS discipline). The waiting space is assumed to be unlimited. Customers are homogeneous with respect to their economic parameters. More concretely, each customer receives a fixed reward $R>0$ upon service completion and incurs a waiting cost at rate $C>0$ per unit of time spent in the system. 
The service provider chooses a nonnegative entrance fee $p<R$. This restriction is without loss of generality, because any fee $p\geq R$ induces all customers to balk and generates zero revenue. We also assume that the service is provided without cost. 

To simplify the analysis and highlight the key trade-offs, we introduce the following dimensionless quantities. 

We define the traffic intensity
\[
\rho = \frac{\lambda}{\mu},
\]
which represents the potential offered load, namely, the traffic intensity that would arise if all potential customers joined. We also define
\[
\nu = \frac{R \mu}{C}, \ \  \nu(p) = \frac{(R - p)\mu}{C}  = \nu - \frac{\mu}{C} p,
\]
which capture the gross and net rewards, respectively, from service relative to the waiting cost. The parameters $\rho$ and $\nu(p)$ provide a convenient and intuitive representation of the system: $\rho$ captures congestion effects, while $\nu$ and $\nu(p)$ reflect the economic attractiveness of joining. This formulation will be used extensively in the analysis that follows.

The customers' dilemma is whether to join or balk. Customers make their decisions without knowing the state of the system (i.e., the model is unobservable), but the type of the system and the various operational and economic parameters of the system, $\lambda,\mu,R, C$ and $p$ are common knowledge.

In this paper we consider two cases regarding the customer decision framework. In the first case, each customer makes the join/balk decision independently. However, since the total arrival rate to the system is affected by the decisions of all customers and in turn affects every customer's waiting cost, this situation gives rise to a game played by all customers, in which symmetric Nash equilibria are sought. We refer to this framework as independent or decentralized behavior.

In the second case the customers are grouped and defer their decisions to a single leader who determines the proportion of them to join the system, with the objective of maximizing the long-run expected aggregate net value of customers per unit time. This framework is referred to as collective or cooperating behavior.

In both cases, given an entrance fee $p$, a fraction $q$ of customers join the system in equilibrium. Then, the system becomes an M/M/1 queue with arrival rate $\lambda q$ and service rate $\mu$; hence the mean sojourn time of a customer who decides to join is equal to 
\[
  W(q) = \frac{1}{\mu - \lambda q} = \frac{1}{\mu (1 - \rho q)},
\]
for $\rho q < 1$  (see e.g., Hassin and Haviv (2003), Section 1.4). 

Depending on the decision framework, the equilibrium joining proportion is denoted by $q_e^{(ind)}(p)$ or $q_e^{(gr)}(p)$.
Therefore, the net value received by a customer who joins the system in equilibrium is equal to $R-p-CW(q_e^{(ind)}(p))$ or $R-p-CW(q_e^{(gr)}(p))$ under independent and grouped behavior, respectively. Note that a balking customer has net value zero. Thus, the customer surplus function, which is defined as the expected total net value of all customers per unit time, is equal to 
\[
  \overline{CS}^{(ind)}(p) = \lambda q_e^{(ind)}(p) \bigl(R - p - C W(q_e^{(ind)}(p))\bigr),
\]
or 
\[
  \overline{CS}^{(gr)}(p) = \lambda q_e^{(gr)}(p) \bigl(R - p - C W(q_e^{(gr)}(p))\bigr).
\]

In addition to the customers' individual or group decisions, we consider the service provider's pricing strategy. The provider acts as a profit maximizer and selects the admission fee $p$ in order to maximize revenue. 

Under each behavioral regime, the effective arrival rate $\lambda q_e(p)$ represents the demand function, i.e., the rate at which customers join the system as a function of the price $p$. Accordingly, the provider's revenue per unit time is given by
\[
\overline{R}^{(ind)}(p) = \lambda q_e^{(ind)}(p) \, p, 
\qquad
\overline{R}^{(gr)}(p) = \lambda q_e^{(gr)}(p) \, p,
\]
under individual and group behavior, respectively.
{Throughout the paper, an overbar denotes a quantity expressed as a function of the entrance fee $p$, whereas the corresponding unbarred quantity is expressed as a function of the induced joining fraction $q$.}

We also consider the total system or social welfare, defined as the aggregate net benefit per unit time of all customers and the service provider. Under any regime, the entrance fee is a transfer payment from the customers to the provider, thus the social welfare depends only on the join probability:
\[
  SW(q) = \lambda q (R - CW(q)) = C \rho q \left ( \nu - \frac{1}{1-\rho q} \right ).
\]

Therefore, in order to maximize the welfare, a social planner should employ a strategy that achieves a value of $q$ that maximizes $SW(q)$. We refer to this value as socially optimal join probability and denote it by $q_{soc}$, i.e.,
\[
  SW(q_{soc}) = SW_{max} \equiv \max_{q \in [0,1]} SW(q) .
\]
Differentiating $SW(q)$ gives the socially optimal joining fraction explicitly:
\begin{equation}
q_{soc}=\min\left\{\frac{1}{\rho}\left(1-\sqrt{\frac{1}{\nu}}\right)^+,1\right\}.
\label{social-optimal-probability}
\end{equation}
When customer behavior is determined in a  strategic equilibrium framework,  the social welfare depends on the entrance fee $p$ indirectly through the equilibrium response $q_e(p)$. This observation allows for two equivalent mechanisms to achieve welfare maximization: either by directly controlling the admission rate (e.g., via access restrictions), or indirectly through pricing.

We define as \emph{coordinating prices} the entrance fees that induce the welfare-maxi-mizing joining proportion, and denote them by $p_{coord}^{(ind)}$ and $p_{coord}^{(gr)}$ under the individual and grouped regimes, respectively.
Specifically  $p_{coord}^{(ind)}$ and $p_{coord}^{(gr)}$ are such that
$$
q_e^{(ind)}(p_{coord}^{(ind)}) = q_e^{(gr)}(p_{coord}^{(gr)}) = q_{soc}.
$$

In the following sections, we analyze each behavioral regime separately. For each case, we characterize the equilibrium joining behavior, the revenue-maximizing price, and the coordinating price. We then compare customer surplus, provider revenue, and total welfare across regimes, and examine how welfare is allocated between customers and the service provider under different pricing objectives.

\section{THE CASE OF INDEPENDENT CUSTOMERS}\label{Section-Anal-Indep-Cust}

In this section, we summarize the findings of Edelson and Hildebrand (1975), who analyzed the independent customer behavior model. The customers are informed about the service fee $p$ set by the service provider, the mean arrival rate $\lambda$ and the mean service rate $\mu$. However they are not informed about the actual number present in the system at the time of their arrival. In this setting, a customer's strategy is generally mixed and is determined by a join probability $q\in[0,1]$. The extreme cases $q=0$ and $q=1$ correspond to the pure strategies balk and join, respectively. 

Since customers are homogeneous and not differentiated in any way, it is reasonable to focus on symmetric equilibrium strategies in the join/balk game. Thus, a join probability $q_e^{(ind)}$ is a symmetric equilibrium under independent customer behavior, if, assuming that all customers follow this strategy, no one has an incentive to deviate from it, or equivalently, $q_e^{(ind)}$ is a best response to itself. 

The equilibrium strategy is determined as follows. Suppose that the population of potential customers follows a join probability $q$. Consider, now, a tagged customer who joins with probability $q'$, when the others join with probability $q$ and the service provider charges entrance fee $p$. Then, her expected utility is 
$$\mathcal{U}^{(ind)}(q';q,p)=(1-q')\cdot 0+ q' \left( R-p-\frac{C}{\mu (1 -\rho q)}\right)=\frac{C}{\mu}\left( \nu(p)-\frac{1}{1-\rho q}\right)q'.$$

Therefore, to find her best response against $q$ and $p$, the tagged customer has to solve the problem $\max_{q'\in[0,1]} \mathcal{U}(q';q,p)$. Since the function $\mathcal{U}(q';q,p)$ is linear with respect to $q'$, the tagged customer bases her decision on the sign of the quantity
$$U^{(ind)}(q,p)=\nu(p)-\frac{1}{1-\rho q}.$$
When $\nu(p)>1$, let
\begin{equation*}
\bar{q}_e^{(ind)}(p)=\frac{1}{\rho}\left( 1-\frac{1}{\nu(p)}\right).
\end{equation*}
be the root of $U^{(ind)}(q,p)$ with respect to $q$. 

We can now proceed to the computation of the equilibrium strategies. First, the strategy `always balk' ($q_e^{(ind)}(p)=0$) is an equilibrium if and only if $U^{(ind)}(0,p)=\nu(p)-1\leq0$, that is, if and only if $\nu(p)\leq1$. Second, a strategy $q_e^{(ind)}(p)\in (0,1)$ is equilibrium strategy, if and only if $q_e^{(ind)}(p) = \bar{q}_e^{(ind)}(p)$, which reduces to $q_e^{(ind)}(p)=\frac{1}{\rho}\left( 1-\frac{1}{\nu(p)}\right)$. This is valid as far as $\bar{q}_e^{(ind)}(p)\in (0,1)$, which occurs if and only if $1< \nu(p)<\frac{1}{1-\rho}$. And, finally, the `always join' ($q_e^{(ind)}(p)=1$) is equilibrium strategy, if and only if $1\leq \bar{q}_e^{(ind)}(p)$, which reduces to $\nu(p)\geq \frac{1}{1-\rho}$. In summary, we have the following result:

\begin{thm}\label{thm-equilibrium-ind}
For the join-or-balk customer dilemma in the unobservable M/M/1 queue with independent customers, a unique equilibrium customer strategy $q_e^{(ind)}(p)$ exists, given by the formula
\begin{eqnarray*}\label{qeind}
q_e^{(ind)}(p)&=&\left\{
\begin{array}{ll}
0, & \nu(p)\leq 1\\
\frac{1}{\rho}\left( 1-\frac{1}{\nu(p)}\right), & 1< \nu(p)<\frac{1}{1-\rho}\\
1, & \nu(p)\geq \frac{1}{1-\rho}
\end{array}
\right.\\
&=&\min\left( \frac{1}{\rho}\left( 1-\frac{1}{\nu(p)}\right)^+,1\right).
\end{eqnarray*}
\end{thm} 
 
Based on Theorem \ref{thm-equilibrium-ind} we can derive the range of achievable input probabilities, as well as the fees that can induce each achievable value.

Note that if $\nu \leq 1$, then $q_e^{(ind)}(p)=0$ for every admissible fee $p\in[0,R)$, since no customer is willing to join even if the system is empty. 

If $\nu >1$, the minimum value of the joining probability is $0$, obtained for $\nu(p) \leq 1$, or equivalently for $p\in\left[\frac{C}{\mu}(\nu - 1),R\right)$.
Since $q_e^{(ind)}(p)$ is decreasing in $p$, its maximum value is achieved for $p=0$ and is equal to 
\begin{equation*}
q_e^{(ind)}(0) = \min \left ( \frac{1}{\rho}\left( 1-\frac{1}{\nu}\right), 1 \right ). 
\end{equation*}

We consider two cases. If $\nu < \frac{1}{1-\rho}$, then $q_e^{(ind)}(0) =  \frac{1}{\rho}\left( 1-\frac{1}{\nu}\right) < 1$, i.e., it is not possible to capture the market and induce all customers to join the system. In this case, 
for any $q\in(0,q_e^{(ind)}(0)]$ the fee that induces a join probability $q$ is unique and  equal to
\begin{equation}\label{ind_price}
p^{(ind)}(q)=\frac{C}{\mu}\left( \nu-\frac{1}{1-\rho q}\right),  
\end{equation}
where without loss of generality we do not consider  fees above $\frac{C}{\mu}(\nu - 1)$ since they all imply zero input rates, revenue and welfare.

On the other hand, if $\nu \geq \frac{1}{1-\rho}$, then  $q_e^{(ind)}(0) = 1$. Then, for $q\in(0,1)$ the fee that induces $q$ is also unique and given by \eqref{ind_price}, whereas $q=1$ can be achieved by any $0\leq p\leq\frac{C}{\mu}\left( \nu-\frac{1}{1-\rho}\right)$. For $q=0$, any $p\in\left[\frac{C}{\mu}(\nu-1),R\right)$ induces zero demand.

We next consider the service provider's revenue maximization problem as well as the customer and social welfare under independent customer behavior.

The service provider by setting the entrance $p$ induces an arrival rate $\lambda_e^{(ind)}(p)=\lambda q_e^{(ind)}(p)$. Conversely, to induce a certain join probability
$q\in(0,q_e^{(ind)}(0)]$, the service provider should charge the entrance fee $p^{(ind)}(q)$ given by \eqref{ind_price}. This is true even in the case $q=1$ where multiple entrance fees achieve market capture, because a revenue-maximizing provider will choose the highest among those. 
Then, the service provider's revenue can be expressed as a function of the induced join probability $q$ as
\begin{equation}\label{revenue-ind-fnct-q}
R^{(ind)}(q)=\lambda q p^{(ind)}(q)=\lambda q \frac{C}{\mu}\left( \nu-\frac{1}{1-\rho q}\right)=C\rho q \left( \nu-\frac{1}{1-\rho q}\right). 
\end{equation}
Note that, since  $q_e^{(ind)}(0) < \frac{1}{\rho}$, $\rho q <1$ for $q \in [0, q_e^{(ind)}(0)]$. Furthermore, 
\begin{eqnarray*}
\frac{d}{dq} R^{(ind)}(q) &=& C\rho \left( \nu - \frac{1}{(1-\rho q)^2} \right),\\
\frac{d^2}{dq^2} R^{(ind)}(q)&=& -\frac{2C\rho^2}{(1-\rho q)^3}<0.
\end{eqnarray*}
Therefore, the function $R^{(ind)}(q)$ is concave for $q\in[0,\frac{1}{\rho})$ and attains its maximum at the root $\bar{q}_{r}^{(ind)}$ of $\frac{d}{dq} R^{(ind)}(q)$ which is given by the formula
\begin{equation}
\bar{q}_{r}^{(ind)}=\frac{1}{\rho}\left( 1-\sqrt{\frac{1}{\nu}}\right).\label{MM1unobservable-social-opt}
\end{equation}
When $\bar{q}_{r}^{(ind)}\in (0,1)$, we deduce that the revenue-maximizing induced joining fraction is given by the formula \eqref{MM1unobservable-social-opt}, otherwise the maximum revenue is attained at $0$ or $1$. More concretely, we have the following result. 

\begin{thm}\label{thm-revenue-max-ind}
For the admission problem in the unobservable M/M/1 queue with independent customers, a unique revenue-maximizing induced joining fraction exists, given by the formula
\begin{eqnarray}
q_{r}^{(ind)}&=& \left\{
\begin{array}{ll}
0, & \nu\leq 1,\\
\frac{1}{\rho}\left( 1-\sqrt{\frac{1}{\nu}}\right), & 1< \nu<\frac{1}{(1-\rho)^2},\\
  1, & \nu\geq \frac{1}{(1-\rho)^2}
\end{array}
\right.
\nonumber \\
&=&\min\left( \frac{1}{\rho}\left(  1-\sqrt{\frac{1}{\nu}} \right)^+,1\right) \label{revenue-max-prob-ind}  
\end{eqnarray}
\end{thm}

The revenue-maximizing joining fraction is induced by the service provider by setting entrance fee
\begin{equation}\label{price-profit-max-ind}
p_r^{(ind)}=p^{(ind)}(q_r^{(ind)})=\left\{
\begin{array}{ll}
0, & \nu\leq 1,\\
\frac{C}{\mu}\left( \nu-\sqrt{\nu}\right), & 1< \nu<\frac{1}{(1-\rho)^2},\\
\frac{C}{\mu}\left( \nu-\frac{1}{1-\rho}\right), & \nu\geq \frac{1}{(1-\rho)^2}.
\end{array}
\right.
\end{equation}
When $\nu\leq1$, the displayed choice $p_r^{(ind)}=0$ is a normalization: every admissible fee $p\in[0,R)$ is revenue maximizing because demand and revenue are zero.
Then, the maximum revenue is given as
\begin{equation}\label{revenue-max-ind}
R_{max}^{(ind)}=\left\{
\begin{array}{ll}
0, & \nu\leq 1,\\
C\left( \sqrt{\nu} -1\right)^2, & 1< \nu<\frac{1}{(1-\rho)^2},\\
C\rho \left( \nu-\frac{1}{1-\rho}\right), & \nu\geq \frac{1}{(1-\rho)^2}.
\end{array}
\right.
\end{equation}
Moreover, we observe that 
$$q_r^{(ind)}\leq q_e^{(ind)}(0),$$ for all parameter values $\rho$, and $\nu$. Indeed, for a mixed strategy $q\in (0,1)$, the zero-fee individual equilibrium condition is $\nu=\frac{1}{1-\rho q}$, whereas the first-order condition for social optimality is $\nu=\frac{1}{(1-\rho q)^2}$. The quadratic term in the latter appears because of the negative externalities of joining (see also Haviv and Oz (2018)). Without an entrance fee, the customers tend to use the system more than what is desirable from the service provider's point of view. 

We next turn to social welfare maximization  and coordination.
From the above analysis it follows that for any achievable join probability $q$ a revenue-maximizing provider sets an entrance fee $p$ such that  $U^{(ind)}(q,p^{(ind)}(q))=0$.
Therefore, the corresponding customers' surplus under revenue maximization is equal to  zero.

We see that, because of the homogeneity of the customers, a revenue-maximizer reaps all the social welfare generated by the system. Therefore, the social welfare, $SW(q)$, is equal to the service provider's revenue, $R^{(ind)}(q)$, so the objectives of a revenue-maximizer and a social-planner coincide:
$$
q_{soc} = q_{r}^{(ind)}
$$
and
$$
SW_{max} = R_{max}^{(ind)}, 
$$
as given in \eqref{revenue-max-prob-ind} and  \eqref{revenue-max-ind}, respectively. 
  
However there is a subtle difference between the two problems regarding the coordinating prices and the social welfare allocation. If $\nu\leq1$, then $q_{soc}=0$ and every admissible fee $p\in[0,R)$ is a coordinating price. Specifically, when $0<q_{soc}<1$, the coordinating price is unique and equal to
$$
p_{coord}^{(ind)} = p_r^{(ind)} = \frac{C}{\mu}\left( \nu-\sqrt{\nu}\right),
$$
thus even under social welfare optimization, the entire welfare is captured by the service provider.  On the other hand, when $q_{soc}=1$, then there is a range of coordinating prices
$$
p_{coord}^{(ind)} \in \left [ 0, \frac{C}{\mu}\left( \nu-\frac{1}{1-\rho}\right) \right ].  
$$
In this case the social optimizer may select any entrance fee $p$ in the above range to impose to the service provider. By doing so, the social welfare
$SW_{max}= C\rho \left( \nu-\frac{1}{1-\rho}\right)$ will be allocated between the service provider and the customers as follows:
$$
\overline{R}^{(ind)}(p) = \lambda p
$$
and
$$
\overline{CS}^{(ind)}(p) = C\rho \left( \nu-\frac{1}{1-\rho}\right) - \lambda p.
$$

\begin{figure}[H]
        \centering
  \caption{Join Probabilities and Coordinating Prices Under Independent Customer Behavior}
  \label{fig:independent} 
  \vspace{16pt}
 \includegraphics[width=\textwidth]{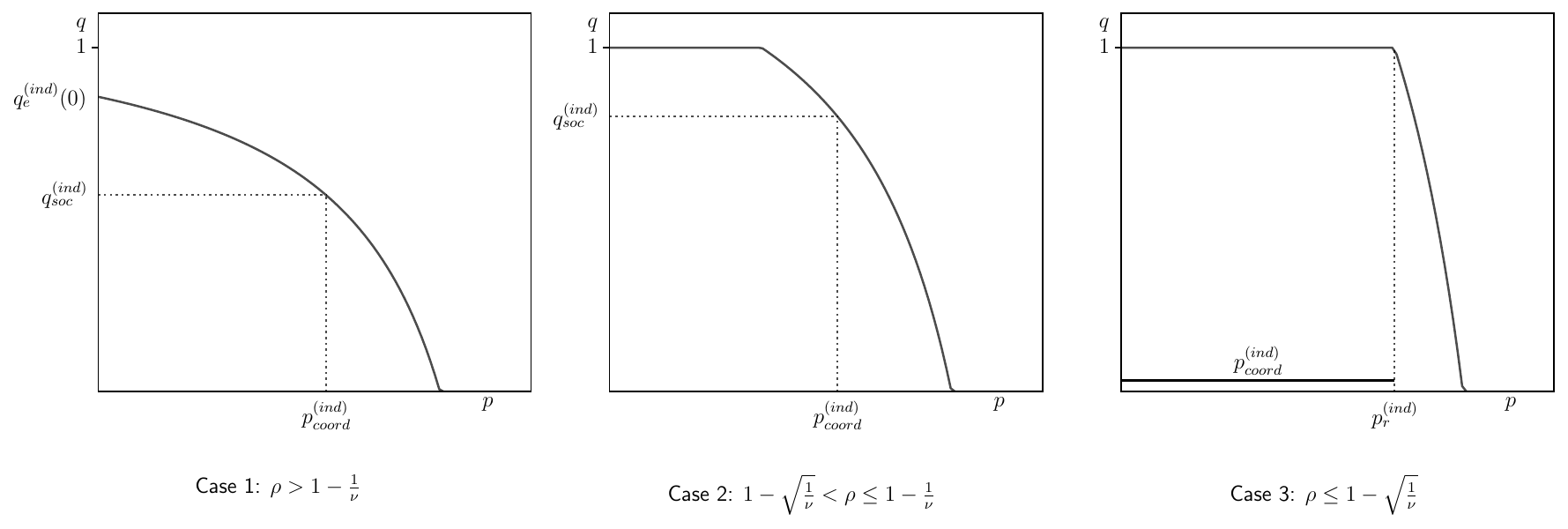}
   \end{figure}

The effect of the traffic intensity  on the join probability and the revenue and social welfare maximization is summarized in  Figure \ref{fig:independent}, which considers three cases for the value of $\rho$. In Case 1, where $\rho > 1- \frac{1}{\nu}$, the customer arrival rate is so high that even with zero entrance fee it is not possible to entice all customers to join. In the intermediate Case 2, where  $1- \frac{1}{\sqrt{\nu}} < \rho \leq 1- \frac{1}{\nu}$, the service provider has the option to capture the market, however it is not profitable to do so, and he still sets the entrance fee at a relatively high level restricting the join probability below 1. In both these cases, the unique coordinating price is the revenue maximizing one which eliminates the customer surplus. This means that a social planner can ensure a positive value for the customers only if he enforces an entrance fee that does not maximize the social welfare. 
Finally in Case 3, where $\rho \leq 1- \frac{1}{\sqrt{\nu}}$ the market size is so low that the optimal solution is to set a market capturing entrance fee. In this case the social planner has significant flexibility in the total welfare allocation between the service provider and the customers since there is a range of coordinating prices.  

\section{THE CASE OF GROUPED CUSTOMERS}\label{Section-Anal-Colluded-Cust}

In the case of grouped customers, we have a classical Stackelberg game, where the first mover is the service provider who sets the entrance fee $p$ and the second mover is the customers' leader who sends a fraction $q$ of them to receive service. The cumulative utility of customers is then 
\begin{equation*}
\mathcal{U}^{(gr)}(q;p)=\lambda q \left( R-p-\frac{C}{\mu-\lambda q} \right)=C\rho q\left(\nu(p)-\frac{1}{1-\rho q}\right).    
\end{equation*}
Unlike $\mathcal{U}^{(ind)}$, which is the expected utility of a tagged customer, $\mathcal{U}^{(gr)}$ is aggregate customer surplus per unit time.
Note that this has essentially the same functional form as \eqref{revenue-ind-fnct-q}, where $\nu$ has been replaced by $\nu(p)$. Therefore, the same analysis is possible with the obvious adaptations. In particular, regarding the best response of the customers' leader, we have the following result.

\begin{thm}\label{thm-equilibrium-col}
The best response of customers' leader for the fraction of potential arrivals that will be sent to the system, in the unobservable M/M/1 queue with grouped customers, is unique and is given by the formula
\begin{eqnarray*}
q_e^{(gr)}(p)&=&\left\{
\begin{array}{ll}
0, & \nu(p)\leq 1\\
\frac{1}{\rho}\left( 1-\sqrt{\frac{1}{\nu(p)}}\right), & 1< \nu(p)<\frac{1}{(1-\rho)^2}\\
1, & \nu(p)\geq \frac{1}{(1-\rho)^2}
\end{array}
\right.\\
&=&\min\left( \frac{1}{\rho}\left( 1-\sqrt{\frac{1}{\nu(p)}}\right)^+,1\right).
\end{eqnarray*}
\end{thm}

Based on the leader's best response function, we can identify the range of achievable join fractions and the corresponding pricing function, similarly to the previous section.

For  $\nu \leq 1$, $q_e^{(gr)}(p)=0$ for every admissible fee $p\in[0,R)$. 
If $\nu >1$, the minimum value of the joining probability is $0$, obtained for $\nu(p) \leq 1$, or equivalently for $p\in\left[\frac{C}{\mu}(\nu - 1),R\right)$.
Since $q_e^{(gr)}(p)$ is decreasing in $p$, its maximum value is achieved for $p=0$ and is equal to 
\begin{equation} \label{eq:qsoc_col}
  q_e^{(gr)}(0) = \min \left ( \frac{1}{\rho}  \left( 1-\sqrt{\frac{1}{\nu} } \right ) , 1  \right ) = q_{soc}. 
\end{equation}
The fact that the join probability under no entrance fee coincides with the socially optimal value has important implications on coordination, which we discuss later in the section.

Regarding  market capture, $q_e^{(gr)}(0)=1$, if and only  if $\nu \geq \frac{1}{(1-\rho)^2}$. Following the same reasoning as in the independent customer regime, we can show that the entrance fee which induces a joining fraction $q \in (0,q_e^{(gr)}(0)]$ is equal to 
\begin{equation}\label{col_price}
p^{(gr)}(q)=\frac{C}{\mu}\left( \nu - \frac{1}{(1-\rho q)^2} \right).
\end{equation}
For $q < 1$ this fee is unique, however if $q_e^{(gr)}(0)=1$, then any $0\leq p \leq \frac{C}{\mu}\left( \nu - \frac{1}{(1-\rho)^2} \right)$ induces market capture. 


For the revenue maximization problem, we express the revenue as a function of the induced join probability $q$, taking into account that when there is more than one choice, the service provider will select the maximum possible price that results in the specific value. Thus, the service provider's objective is the maximization of 
\begin{equation}\label{revenue-col-fnct-q}
R^{(gr)}(q)=\lambda q p^{(gr)}(q)=\lambda q \frac{C}{\mu}\left( \nu-\frac{1}{(1-\rho q)^2}\right)=C\rho q \left( \nu-\frac{1}{(1-\rho q)^2}\right).
\end{equation}
 We have that
\begin{eqnarray*}
\frac{d}{dq} R^{(gr)}(q) &=& C\rho \left( \nu - \frac{1+\rho q}{(1-\rho q)^3} \right),\\
\frac{d^2}{dq^2} R^{(gr)}(q)&=& -\frac{2C\rho^2(2+\rho q)}{(1-\rho q)^4}<0.
\end{eqnarray*}
We conclude that the function $R^{(gr)}(q)$ is concave for $q\in[0,\frac{1}{\rho})$. The first-order condition $\frac{d}{dq} R^{(gr)}(q)=0$ is equivalent to  
$$\nu=\frac{1+\rho q}{(1-\rho q)^3}.$$
Setting $y=1-\rho q$, it can be written as
$$\nu y^3+y-2=0.$$
Define
\begin{equation*}
h(y) = \nu y^3 + y - 2.
\end{equation*}
We have that $h'(y) = 3\nu y^2 + 1 > 0$, so $h(y)$ is strictly increasing and thus has at most one real root. Moreover, $h(0) = -2 < 0$, $h(1) = \nu - 1$, so we have to consider two cases according to whether $\nu > 1$ or $\nu \leq 1$.

If $\nu > 1$, then $h(1) > 0$, hence there is exactly one root of $h(y)$ in $(0,1)$. In the complementary case where $\nu \leq 1$, we have that $h(1)\leq 0$, so there is no root of $h(y)$ in $(0,1)$.

The root can be given explicitly. To this end, we divide the equation by $\nu$ to transform it to a depressed cubic with unit leading coefficient, i.e. to $y^3 + a y + b = 0$, with $a = \frac{1}{\nu}$ and $b = -\frac{2}{\nu}$.
Using the cubic formula
\begin{equation*}
y = \sqrt[3]{-\frac{b}{2} + \sqrt{\left(\frac{b}{2}\right)^2 + \left(\frac{a}{3}\right)^3}}
+ \sqrt[3]{-\frac{b}{2} - \sqrt{\left(\frac{b}{2}\right)^2 + \left(\frac{a}{3}\right)^3}}.
\end{equation*}
we have that $y=y(\nu)$  can be given explicitly and we have the following result. The cube roots in this expression are understood as real cube roots.

\begin{thm}\label{thm-revenue-max-col}
For the admission problem in the unobservable M/M/1 queue with grouped customers, the provider's unique revenue-maximizing induced joining fraction is
\begin{eqnarray*}
q_{r}^{(gr)}&=&\left\{
\begin{array}{ll}
0, & \nu\leq 1,\\
\frac{1}{\rho}\left( 1-y(\nu)\right), & 1<\nu<\frac{1+\rho}{(1-\rho)^3},\\
1, & \nu\geq \frac{1+\rho}{(1-\rho)^3} 
\end{array}
\right. \\
&=& \min \left ( \frac{1}{\rho}\left( 1-y(\nu) \right )^+ , 1 \right ),
\end{eqnarray*}
where
\begin{equation}
y(\nu) =
\sqrt[3]{\frac{1}{\nu} + \sqrt{\frac{1}{\nu^2} + \frac{1}{27\nu^3}}}
+
\sqrt[3]{\frac{1}{\nu} - \sqrt{\frac{1}{\nu^2} + \frac{1}{27\nu^3}}}.
\end{equation}
\end{thm}

The revenue-maximizing joining fraction is imposed by the service provider by setting the entrance fee to 
\begin{equation}\label{price-profit-max-col}
p_r^{(gr)}=p^{(gr)}(q_r^{(gr)})=\left\{
\begin{array}{ll}
0, & \nu\leq 1,\\
\frac{C}{\mu}\left( \nu-\frac{1}{y(\nu)^2}\right), & 1<\nu< \frac{1+\rho}{(1-\rho)^3},\\
\frac{C}{\mu}\left( \nu-\frac{1}{(1-\rho)^2}\right), & \nu\geq \frac{1+\rho}{(1-\rho)^3}.
\end{array}
\right.
\end{equation}
When $\nu\leq1$, the displayed choice $p_r^{(gr)}=0$ is a normalization: every admissible fee $p\in[0,R)$ is revenue maximizing because demand and revenue are zero.
Then the maximum revenue is given by
\begin{equation}\label{revenue-max-col}
R_{max}^{(gr)}=\left\{
\begin{array}{ll}
0, & \nu\leq 1,\\
C(1-y(\nu))\left( \nu-\frac{1}{y(\nu)^2}\right), & 1<\nu < \frac{1+\rho}{(1-\rho)^3},\\
C\rho \left( \nu-\frac{1}{(1-\rho)^2}\right), & \nu\geq \frac{1+\rho}{(1-\rho)^3}.
\end{array}
\right.
\end{equation}

The customers' surplus under the maximal entrance fee that induces a given joining fraction $q$ is 
\begin{equation*}
  CS^{(gr)}(q)=\mathcal{U}^{(gr)}(q;p^{(gr)}(q))=C\rho^2 \frac{q^2}{(1-\rho q)^2} =
  C \left ( \frac{1}{1-\rho q} -1 \right )^2.
\end{equation*}
Therefore, in the case of grouped customers, a revenue-maximizing provider cannot reap all the social welfare as in the independent customers case.
Note that, as expected, the social welfare function is the same as in the independent regime:
\begin{equation*}
  SW^{(gr)}(q) = R^{(gr)}(q) + CS^{(gr)}(q) = C \rho q \left ( \nu -\frac{1}{1-\rho q} \right ).
\end{equation*}
It was shown in \eqref{eq:qsoc_col} that $q_{soc}=q_e^{(gr)}(0)$; thus, the social planner is now faced with the opposite problem compared to the independent behavior case. If $\nu\leq1$, then $q_{soc}=0$ and every admissible fee $p\in[0,R)$ is a coordinating price. Specifically, if $0<q_{soc}<1$, then the unique coordinating price is $p_{coord}^{(gr)}=0$, so the only way to maximize the total welfare is to force the provider to set the entrance fee to zero, thus eliminating his own revenue. If $q_{soc}=1$, then the social planner has more flexibility; any fee that induces all customers to join leads to coordination:
\begin{equation*}
p_{coord}^{(gr)}\in \left[0,\frac{C}{\mu}\left(\nu-\frac{1}{(1-\rho)^2}\right)\right].
\end{equation*}

If the service provider is allowed to charge the revenue maximizing fee, the welfare of the two  parties is determined as follows.
First, for $\nu \leq 1$, no customers join the system, thus $R_{max}^{(gr)} = CS^{(gr)}(q_r^{(gr)}) = 0$. When
$1<\nu < \frac{1}{(1-\rho)^2}$, the provider's revenue is 
\[
R_{max}^{(gr)} =   C(1-y(\nu))\left( \nu-\frac{1}{y(\nu)^2}\right)
\]
and the maximum social welfare
\[
SW_{max} = C\left( \sqrt{\nu} -1\right)^2.
\]
Letting $\rho q_r^{(gr)} = 1 - y(\nu)$ in the customer surplus function we obtain
\begin{equation}
  \label{eq:customer_surplus_col}
  CS^{(gr)}(q_r^{(gr)}) = C \frac{\left (\rho q_r^{(gr)}\right )^2}{(1-\rho q_r^{(gr)})^2}
  = C \left ( \frac{1}{y(\nu)} -1 \right )^2,  
\end{equation}

Therefore, the social welfare under revenue maximization is equal to
\[
  SW(q_r^{(gr)}) = C \rho q_r^{(gr)} \left ( \nu -\frac{1}{1-\rho q_r^{(gr)}} \right )
  = C(1-y(\nu)) \left (\nu - \frac{1}{y(\nu)} \right ).
\]
The relative loss in social welfare under revenue maximization is defined, for $\nu>1$, as
\begin{equation*}
L^{(gr)} =   \frac{SW_{max} -  SW(q_r^{(gr)})}{SW_{max}}.
\end{equation*}
Using the expression for $SW_{max}$ and the preceding expression for $SW(q_r^{(gr)})$, we obtain
\begin{equation*}
  L^{(gr)} =   \frac{(\sqrt{\nu} - 1)^2 - (1-y(\nu)) \left (\nu - \frac{1}{y(\nu)} \right )}{(\sqrt{\nu} - 1)^2}
  =  \frac{  \nu y(\nu) + \frac{1}{y(\nu)}  - 2\sqrt{\nu} }{(\sqrt{\nu} - 1)^2}
\end{equation*}
However, from the cubic equation that $y(\nu)$ solves we obtain:
\[
   \nu y(\nu) + \frac{1}{y(\nu)} = \frac{1}{y(\nu)^2} (\nu y(\nu)^3 + y(\nu)) = \frac{2}{y(\nu)^2}, 
\]
thus,
\begin{equation*}
  L^{(gr)} =   \frac{2 \left (  \frac{1}{y(\nu)^2}  - \sqrt{\nu} \right ) }{(\sqrt{\nu} - 1)^2}. 
\end{equation*}

When $ \frac{1}{(1-\rho)^2} \leq \nu < \frac{1+\rho}{(1-\rho)^3}$, we still have
\[
  SW^{(gr)}(q_r^{(gr)}) =  C(1-y(\nu)) \left (\nu - \frac{1}{y(\nu)} \right ),
\]
however the maximum social welfare now is equal to
\[
  SW_{max} = C \rho \left (\nu - \frac{1}{1-\rho} \right ),
\]
therefore the relative welfare loss becomes
\begin{equation*}
  L^{(gr)} =   1 - \frac{ (1-y(\nu)) \left (\nu - \frac{1}{y(\nu)} \right ) }{\rho \left (\nu - \frac{1}{1-\rho} \right )}. 
\end{equation*}

Finally, when $\nu \geq \frac{1+\rho}{(1-\rho)^3}$, we have  $q_r^{(gr)}=1$ and as it was shown above the revenue maximizing fee is coordinating, thus $SW^{(gr)}(q_r^{(gr)}) = SW_{max}$ and $L^{(gr)}=0$.

Summarizing, the relative loss in social welfare when the provider maximizes revenue is equal to
\begin{equation}
 \label{eq:welfare_loss_col}
  L^{(gr)} = \left \{
\begin{array}{ll}
\frac{2 \left (  \frac{1}{y(\nu)^2}  - \sqrt{\nu} \right ) }{(\sqrt{\nu} - 1)^2}, &  1< \nu<\frac{1}{(1-\rho)^2}\\
  1 - \frac{ (1-y(\nu)) \left (\nu - \frac{1}{y(\nu)} \right ) }{\rho \left (\nu - \frac{1}{1-\rho} \right )} & \frac{1}{(1-\rho)^2} \leq \nu < \frac{1+\rho}{(1-\rho)^3} \\
0, &\nu \geq \frac{1+\rho}{(1-\rho)^3}
\end{array}
\right.\\  
\end{equation}

\begin{figure}[H] 
        \centering
  \caption{Join Probabilities and Coordinating Prices Under Grouped Customer Behavior}
  \label{fig:colluded} 
  \vspace{16pt}
 \includegraphics[width=\textwidth]{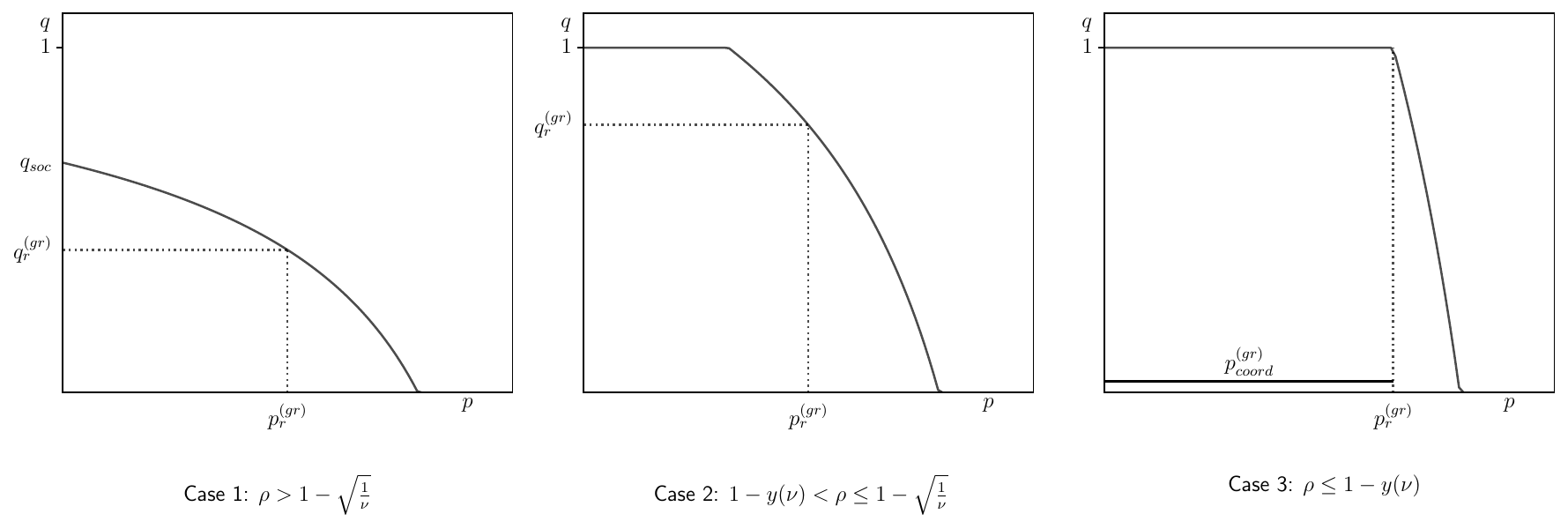}
\end{figure}

Figure \ref{fig:colluded} summarizes the effect of the traffic intensity on the revenue maximizing strategy and the coordinating prices. In Cases 1 and 2 the potential offered load is so large that the provider is either not able or not willing to make all customers join. In addition, in both cases the revenue maximization strategy is suboptimal with respect to the social welfare. In Case 3 where the potential offered load is low, both the revenue maximizing and the socially optimal entrance fee are market capturing, and the system can be coordinated without eliminating either party's surplus. 

\section{INDEPENDENT VERSUS GROUPED CUSTOMERS}\label{Section-Comparison}

In this section we compare the cases of  independent and grouped customer behavior, in terms of the equilibrium join probabilities, revenue-maximizing and coordinating fees and the welfare of both parties.

We first compare the equilibrium join probabilities derived in Theorems \ref{thm-equilibrium-ind} and \ref{thm-equilibrium-col}. When $\nu(p) \leq 1$ both $q_e^{(ind)}(p) = q_e^{(gr)}(p) = 0$. When $\nu(p)>1$, we have $\frac{1}{\nu(p)} < \sqrt{\frac{1}{\nu(p)}}$, thus, $q_e^{(ind)}(p) \geq q_e^{(gr)}(p)$.
Therefore, a larger fraction of customers join the system under independent than under group behavior. This is expected, since in the independent case each customer considers only her own benefit from joining the system, whereas in the grouped case the customer leader also takes into account the negative externalities, in terms of the extra delay that each incoming customer induces to future arrivals.

Regarding the revenue-maximizing probabilities $q_r^{(ind)}$ and $q_r^{(gr)}$, for $\nu\leq 1$ both probabilities are $0$. Suppose that $\nu>1$. Then $y(\nu)<1$ (since it is the unique solution of $\nu y^3+y-2=0$ in $(0,1)$). Therefore,
\begin{eqnarray*}
&&2 y(\nu)<2 \Rightarrow 2 y(\nu)<\nu y(\nu)^3+y(\nu)\\
&&\Rightarrow y(\nu)<\nu y(\nu)^3\Rightarrow 1<\nu y(\nu)^2 \Rightarrow \sqrt{\frac{1}{\nu}}<y(\nu).
\end{eqnarray*}
Hence $$q_r^{(ind)}=\min\left(\frac{1}{\rho}\left( 1-\sqrt{\frac{1}{\nu}}\right)^+,1\right)\geq \min\left(\frac{1}{\rho}\left( 1-y(\nu)\right)^+,1\right)=q_r^{(gr)}. $$
Accordingly, even when the entrance fee is set to maximize revenue, grouped behavior leads to a weakly lower joining fraction.

We will now compare the revenue-maximizing prices. For $\nu\leq 1$, every admissible price $p\in[0,R)$ is revenue maximizing in both regimes; we select $0$ as a normalization. Thus, we focus on the case where $\nu>1$. 

Note that we can write \eqref{price-profit-max-ind} as
\begin{equation*}
p_r^{(ind)}=\left\{
\begin{array}{ll}
\frac{C}{\mu}\left(\nu-\frac{1}{1-\rho}\right), & \mbox{ if } \rho<1-\sqrt{\frac{1}{\nu}}\\
\frac{C}{\mu}\left(\nu-\sqrt{\nu}\right), & \mbox{ if } \rho\geq 1-\sqrt{\frac{1}{\nu}},
\end{array}
\right.
\end{equation*}
which is a decreasing and ultimately constant function of $\rho$. This may seem a paradox at first glance, since $\rho$ can be interpreted as the relative demand for service with respect to the service capacity,  therefore we might expect that an increase in the demand would make the monopolist raise the price. However, an increase in the demand induces a significant decrease in the `quality' of the service, because of increasing delays. Therefore, the customers become more reluctant to buy the service and the monopolist cannot increase the price.

A similar situation occurs in the grouped customers' case. Then \eqref{price-profit-max-col} can be written as

\begin{equation*}
p_r^{(gr)}=\left\{
\begin{array}{ll}
\frac{C}{\mu}\left( \nu-\frac{1}{(1-\rho)^2}\right), & \rho<1-y(\nu)\\
\frac{C}{\mu}\left( \nu-\frac{1}{y(\nu)^2}\right), & \rho\geq 1-y(\nu).
\end{array}
\right.
\end{equation*}
By considering the cases $\rho\leq 1-y(\nu)$, $1-y(\nu)<\rho<1-\sqrt{\frac{1}{\nu}}$ and $\rho\geq 1-\sqrt{\frac{1}{\nu}}$ and comparing the corresponding branches of $p_r^{(ind)}$ and $p_r^{(gr)}$, we can easily verify that in all cases
$$p_r^{(gr)}\leq p_r^{(ind)}.$$

We see that, in order to maximize his revenue, the service provider is forced to set a lower price in the grouped customers' case and, even then, fewer customers do enter the system. It is then  obvious that customer cooperation is generally detrimental for the service provider's revenue. An interesting question is the behavior of the relative loss in revenue due to customers' cooperation with respect to the parameters $\rho$ and $\nu$, namely, for $\nu>1$, the ratio
$$L_r = \frac{R_{max}^{(ind)}-R_{max}^{(gr)}}{R_{max}^{(ind)}}.$$
Recall that $R_{max}^{(ind)} = SW_{max}$ and $R_{max}^{(gr)} = SW(q_r^{(gr)}) - CS^{(gr)}(q_r^{(gr)})$. Therefore,
\begin{equation*}
  L_r = L^{(gr)} + \frac{ CS^{(gr)}(q_r^{(gr)})}{SW_{max}}, 
\end{equation*}
where $L^{(gr)}$ is the loss in social welfare in the grouped regime defined in \eqref{eq:welfare_loss_col}, from which we obtain
\begin{equation}
 \label{eq:revenue_loss_col}
  L_r = \left \{
\begin{array}{ll}
  \frac{2 \left (  \frac{1}{y(\nu)^2}  - \sqrt{\nu} \right ) + \left ( \frac{1}{y(\nu)} -1 \right )^2 }{(\sqrt{\nu} - 1)^2},
  &  1< \nu<\frac{1}{(1-\rho)^2}\\

  1 - \frac{ (1-y(\nu)) \left (\nu - \frac{1}{y(\nu)} \right ) - \left ( \frac{1}{y(\nu)} -1 \right )^2 }
  {\rho \left (\nu - \frac{1}{1-\rho} \right )},
  & \frac{1}{(1-\rho)^2} \leq \nu < \frac{1+\rho}{(1-\rho)^3} \\

  \frac{\rho}{(1-\rho) ( \nu (1-\rho) - 1)}
  ,
  &\nu \geq \frac{1+\rho}{(1-\rho)^3}
\end{array}
\right.\\  
\end{equation}

We see that even in the case $\nu \geq \frac{1+\rho}{(1-\rho)^3}$ where the revenue maximizing fee is such that the market is captured and the system is coordinated, the customers still receive a positive part of the social welfare.

Finally, regarding the comparison of coordinating prices, it follows from the previous sections that when the socially optimal join probability satisfies $0<q_{soc}<1$, then $p_{coord}^{(gr)} = 0 < p_{coord}^{(ind)}$. When $q_{soc}=0$, every admissible price $p\in[0,R)$ coordinates both regimes. When $q_{soc}=1$, then the service provider can induce coordination in both regimes by setting the entrance fee at or below a maximum value. In this case the maximum value to achieve coordination is lower under grouped than under independent behavior. 

\section{COMPUTATIONAL RESULTS}\label{Section-Computational}

In this section we discuss the sensitivity of  measures of performance discussed above with respect to the two key parameters of the problem, the traffic intensity $\rho$ and the relative service reward $\nu$. In the computational experiments we have set the waiting cost per unit time $C$ and the service rate $\mu$ to nominal values of 1, and varied $\rho$ in the range $(0,1)$ and $\nu$ in the range $[2,20]$.

\begin{figure}[H]
        \centering
        \caption{Comparison of Join Probabilities, Fees and Revenues Between the Two Regimes for $\nu=2$}
  \label{fig:comp_nu_2} 
  \vspace{16pt}
 \includegraphics[width=\textwidth]{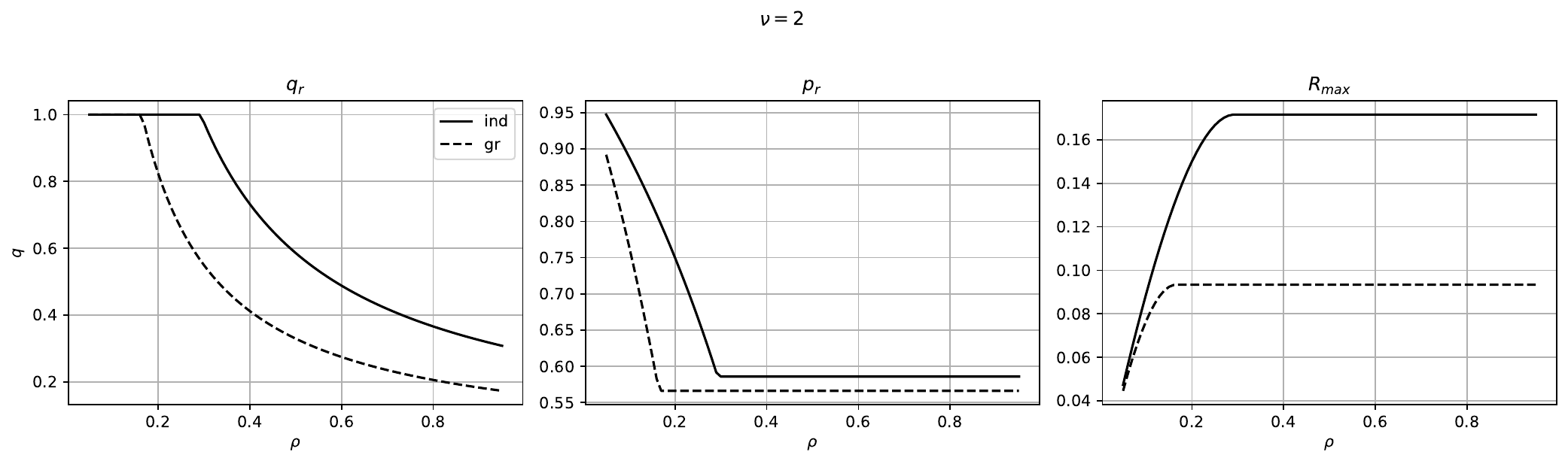}
\end{figure}

\begin{figure}[H]
        \centering
        \caption{Comparison of Join Probabilities, Fees and Revenues Between the Two Regimes for $\nu=10$}
  \label{fig:comp_nu_10} 
  \vspace{16pt}
 \includegraphics[width=\textwidth]{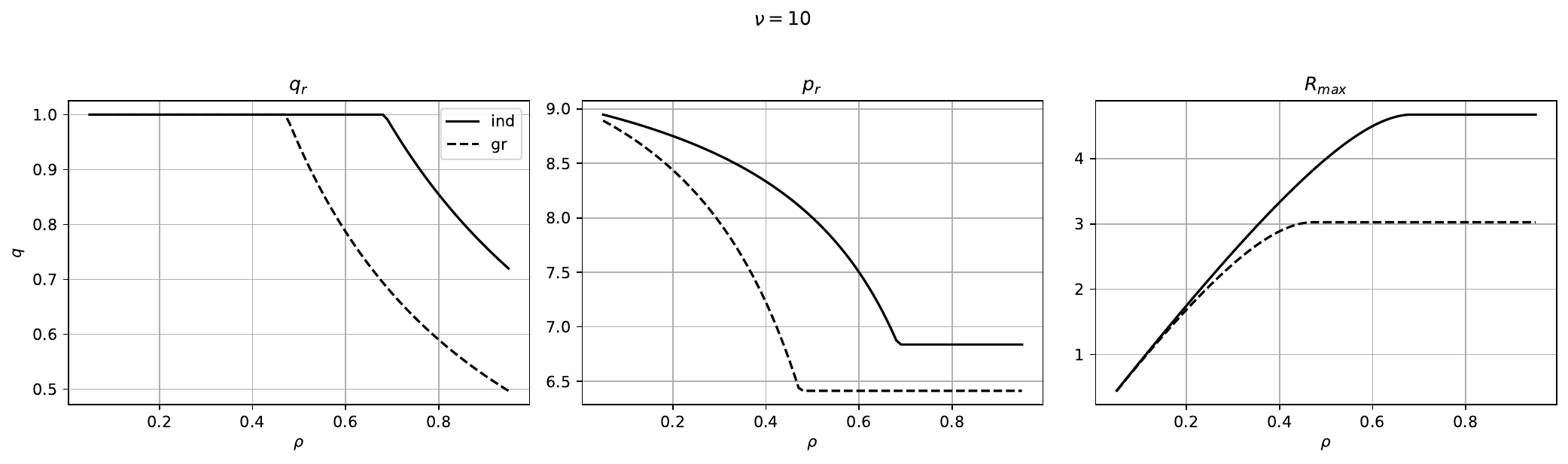}
\end{figure}

In Figures \ref{fig:comp_nu_2} and  \ref{fig:comp_nu_10} we present the comparison between the independent and grouped customer regimes, in terms of the revenue maximizing join probabilities $q_r$, service fees $p_r$ and maximum revenues $R_{max}$ as a function of the traffic intensity, for  $\nu =2$ and $\nu=10$, respectively.

Several observations can be made here. First, for small values of the traffic intensity the market is captured ($q_r = 1$), although as $\rho$ increases, the price required to keep all customers joining falls. Nevertheless the price drop is relatively lower than the increase in demand for service, thus the maximum revenue also increases. After the traffic-intensity threshold at which the market is no longer captured, the joining fraction $q_r$ decreases as $1/\rho$, while the effective load $\rho q_r$, the revenue-maximizing fee, and the maximum revenue become constant.

Increasing the relative service value $\nu$ raises the revenue-maximizing joining fraction, fee, and revenue in both regimes.

Regarding the comparison between the two regimes, the plots verify the observations in the previous sections. When customers are able (or forced) to coordinate in terms of their response to the provider's pricing strategy, this cooperation is detrimental for the service provider's revenue. Indeed, the provider is forced to lower the entrance fee significantly, and even then the arrival rates in equilibrium are lower, which is reflected in lower revenues.

\begin{figure}[H]
        \centering
        \caption{Welfare Loss and Revenue Loss Under Grouped Regime}
  \label{fig:loss} 
  \vspace{16pt}
 \includegraphics[width=\textwidth]{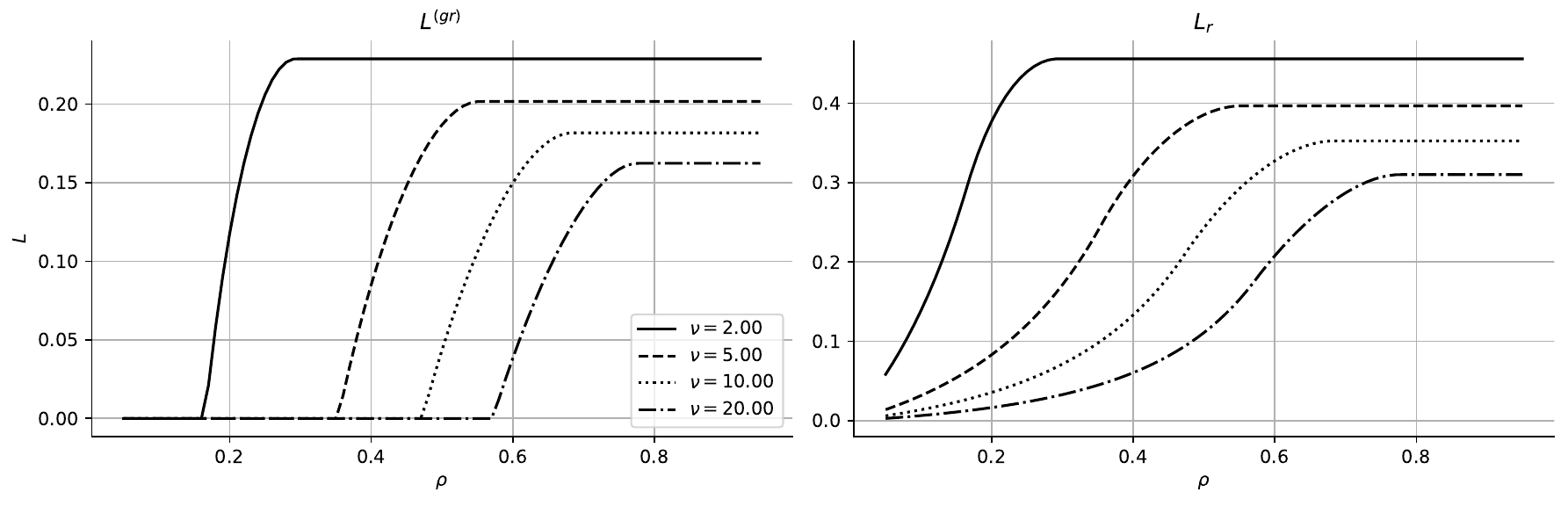}
\end{figure}

Figure \ref{fig:loss} presents the relative welfare loss $L^{(gr)}$ and the relative revenue loss between grouped and independent behavior, $L_r$, as functions of $\rho$ for several values of $\nu$.

The left panel presents the sensitivity of the relative welfare loss under cooperation $L^{(gr)}$. 
As long as the market is captured under revenue maximization the loss is zero since the social welfare is maximized. When $\rho$ increases beyond the threshold at which the service provider no longer prefers to capture the market but still the socially optimal probability $q_{soc} = q_r^{(ind)} =1$,  the relative welfare loss increases steeply. When $q_{soc}$ falls below 1 the loss becomes constant since the welfare values are stabilized. When the service value increases, the relative loss is reduced but not significantly compared to the increase in $\nu$.

The right panel presents the sensitivity of the relative revenue reduction between independent and grouped regimes,  $L_r$. The most important difference from the previous panel is that here the loss in the provider's revenue is still positive for small values of $\rho$. This is so because, when customers behave cooperatively, even when the social welfare is maximized under market capture, the service provider is not able to reap all customer surplus and his revenue is still lower than under independent customer behavior.  

\section{CONCLUSIONS}\label{Section-Conclusions}

In this paper we examined the join--or--balk decision in an unobservable single--server queue under two alternative behavioral regimes: independent customers and grouped  (cooperating) customers. The comparison highlights how the nature of decision making fundamentally affects equilibrium participation, system congestion, and pricing outcomes.

The analysis shows that cooperation among customers leads to systematically different behavior. In contrast to independent decision making, where each customer ignores the congestion imposed on others, coordinated customers internalize this externality and therefore adopt more conservative joining strategies. As a result, cooperating behavior leads to a weakly lower joining fraction and effective arrival rate at any given fee.

These differences propagate to the provider’s pricing problem. Since demand is endogenously determined by customer behavior, the optimal pricing policy depends critically on whether customers act independently or in a coordinated manner. In particular, coordination alters both the shape of the demand function and the resulting revenue-maximizing price. The revenue-maximizing joining fraction, fee, and provider revenue are all weakly lower under grouped behavior.

The welfare comparison is more specific. Under independent behavior, the provider's revenue-maximizing joining fraction coincides with the socially optimal fraction, and the provider extracts all surplus except when a social planner selects a lower coordinating fee in the market-capture region. Under grouped behavior, a zero fee implements the social optimum, whereas the provider's revenue-maximizing fee generally induces a joining fraction below the social optimum and leaves customers with positive surplus. For $\nu>1$, social efficiency under grouped revenue maximization is recovered only when the revenue-maximizing policy captures the market; for $\nu\leq1$, both the revenue-maximizing and socially optimal joining fractions are zero.

Overall, the results demonstrate that incorporating customer cooperation into queueing models provides a richer and more nuanced understanding of strategic behavior in service systems. Even in the simple setting considered here, the contrast between independent and group decisions leads to qualitatively different outcomes, underscoring the importance of explicitly modeling the decision structure of customers.

Several directions for future research arise naturally from the present analysis. An important extension is the study of partially observable systems, where the interaction between information and cooperation may lead to richer equilibrium structures. It would also be of interest to examine intermediate forms of cooperation, where customers are only partially aligned, as well as dynamic models in which cooperation evolves over time.

\section{FIGURES}
 
\begin{center}
{\large {\bf REFERENCES} }
\end{center}
\vspace*{4pt}
\begin{description}
\item Benioudakis, M., Zissis, D., Burnetas, A. and Ioannou, G. (2023) Service provision on an aggregator platform with time-sensitive customers: Pricing strategies and coordination,  \textit{International Journal of Production Economics}, \textbf{257}, 108760.

\item Bountali, O., Burnetas, A. and \"{O}rmeci, L. (2022)
  Join, balk, or jettison? The effect of flexibility and ranking knowledge in systems with batch arrivals\textcolor{red}{,}
  \textit{Production and Operations Management},
  \textbf{31}, 3505-3524.

\item Economou, A. (2020) The impact of information structure on strategic behavior in queueing systems, in Anisimov, V. and Limnios, N. \textit{Queueing Theory 2: Advanced Trends}, ISTE Ltd and John Wiley and Sons, Inc. Chapter 4, 137--168.

\item Edelson, N.M. and Hildebrand, K. (1975) Congestion tolls for Poisson queueing processes, \textit{Econometrica} \textbf{43}, 81--92.

\item Hassin, R. (2016) \textit{Rational Queueing}, CRC Press, Taylor and Francis Group, Boca Raton.

\item Hassin, R. and Haviv, M. (2003) \textit{To Queue or Not to Queue: Equilibrium Behavior in Queueing Systems}, Kluwer Academic Publishers, Boston.

\item Haviv, M. and Oz, B. (2018) Self-regulation of an unobservable queue, \textit{Management Science}, \textbf{64}(5), 2380--2389.

\item Naor, P. (1969) The regulation of queue size by levying tolls, \textit{Econometrica}, \textbf{37}, 15--24.

\item Stidham, S. Jr. (2009) \textit{Optimal Design of Queueing Systems}, CRC Press, Taylor and Francis Group, Boca Raton.

\item Tirole, J. (1988) \textit{The Theory of Industrial Organization},  MIT Press, Cambridge, MA.

\end{description}
        
\end{document}